\documentclass[12pt]{amsart}
\usepackage{amssymb,amsmath,amsthm}
\def\pf{ \noindent {\bf Proof: \  }}
\renewcommand{\qed}{\hfill\vrule height6pt  width6pt depth0pt}

\font\Bbb=msbm10 scaled \magstep1 \font\scriptBbb=msbm10
\font\scriptscriptBbb=msbm7
\newfam\Bbbfam
\textfont\Bbbfam=\Bbb \scriptfont\Bbbfam=\scriptBbb
\scriptscriptfont\Bbbfam=\scriptscriptBbb
\def\beq{\begin{equation}}

\def\tr{\mathrm{tr\ }}

\def\calf{{\mathcal{F}}}
\def\call{{\mathcal{L}}}

\newtheorem{theorem}{Theorem}[section]
\newtheorem{lemma}{Lemma}
\newtheorem{prop}{Proposition}
\newtheorem{prop*}{Proposition}
\newtheorem{corollary}{Corollary}
\theoremstyle{definition}

\theoremstyle{remark}

\def\a{\alpha}

\begin{document}
\title[Dual Form of  AP]{The Dual Form of the Approximation Property for a 
Banach Space and a Subspace}\author{T.~Figiel}
\address{Institute of Mathematics, The Polish Academy of Sciences} 
\email{t.figiel@impan.gda.pl} \author{W.~B.
~Johnson$^*$} \address{Department of Mathematics, Texas A\&M
University, College Station, TX 77843--3368 U.S.A}
\email{johnson@math.tamu.edu}
 \date{August 4,
2015.}
\thanks{$^*$ Supported in part by NSF
DMS-1310550}

\begin{abstract}
Given a Banach space $X$ and a subspace $Y$, the pair $(X,Y)$ is said to have the approximation
property (AP) provided there is a net of 
finite rank bounded linear operators on $X$ all of which leave the 
subspace $Y$ invariant such that the net converges uniformly on compact subsets of $X$ to the identity operator.
 The main result is an easy to apply dual formulation of this property. Applications are given to three space properties; in particular, if $X$ has the approximation property and its subspace $Y$ is  $\call_\infty$, then $X/Y$ has the approximation property.
\end{abstract}

\maketitle

  \begin{center}
In memory of A. Pe\l czy\'nski
\end{center}

\setcounter{section}{0}

\section{Introduction}\label{introduction}

In \cite{FJP} the authors and the late A. Pe\l czy\'nski introduced the 
notion of the bounded approximation property (BAP) for a Banach space $X$ 
and a subspace $Y$.  The pair $(X,Y)$ is said to have the approximation 
property (AP) provided the identity on $X$ is the $\tau$-limit of a net of 
finite rank bounded linear operators on $X$ all of which leave the 
subspace $Y$ invariant.  Here we recall that the $\tau$-topology on the 
space $L(X)$ of bounded linear operators on the Banach space $X$  is the 
topology of uniform convergence on compact subsets of $X$.  If the 
approximating net of finite rank operators can be chosen so that their 
norms are uniformly bounded by $\lambda$, then $(X,Y)$ is said to have the 
$\lambda$-bounded approximation property ($\lambda$-BAP), and $(X,Y)$ has 
the BAP provided it has the $\lambda$-BAP for some $\lambda<\infty$. When 
the subspace $Y$ is either the whole space or the zero subspace, these 
concepts reduce to the classical concepts of AP and BAP for a single 
space.  Obviously if $(X,Y)$ has the BAP then $(X,Y)$ has the AP.  Non 
obvious is the fact, pointed out by Lissitsin and Oja \cite[Corollary 
5.12]{LO}, that if $X$ is reflexive and $(X,Y)$ has the AP, then $(X,Y)$ 
has the $1$-BAP.  The root for this is Grothendieck's classical theorem 
\cite{Gr}, that the AP implies the $1$-BAP for reflexive spaces,  as 
improved by
Godefroy-Saphar \cite[Theorem 1.5]{GS}, Oja and collaborators (see, e.g., 
\cite{LMO} and references therein), and others.
  The AP for a pair $(X,Y)$ was not considered in \cite{FJP}, but even if 
one cares only about reflexive spaces it is probably worthwhile to 
consider the concept, in part because the dual form of the AP is simply 
stated and easy to work with while the dual form of the BAP is more 
problematical.   In this note we prove in Theorem \ref{JointAPthm}    the 
dual form for what it means for a pair $(X,Y)$ to have the AP and give a 
couple of applications.  In a paper under preparation we will give a far 
reaching extension of the duality result where the subspace $Y$ is 
replaced by a nest of subspaces of $X$.  Although the proof of the general 
result is not essentially more complicated than what is treated here, it 
does require introducing concepts extraneous to the context of this short note. 
It seemed to us that the special case considered here as well as the 
applications were   interesting enough to warrant a separate publication. The 
applications of  Theorem \ref{JointAPthm} are new approximation property 
three space results in the spirit of other such results (see e.g. 
\cite{GS2} and  \cite{CK}).

We use standard Banach space theory notation and concepts, as are 
contained e.g. in \cite{LT}.

\section{Joint AP}\label{JointAP}

We begin with a special case of a known lemma (part (1) is contained
 in \cite{ring2} and part (2) is in 
\cite{Sp}) but include a simple proof.

\begin{lemma}\label{lemmaRank1}
Let $\calf_Y(X) = \{ T \in \calf(X) : TY \subseteq Y \}$.
\item{1.}  $x^*\otimes x \in \calf_Y(X) $ if and only if either $x^*\in 
Y^\perp$ or $ x\in Y$.
\item{2.} If  $F\in \calf_Y(X) $, then $F$  is the sum of  $n$ rank one elements of 
$\calf_Y(X)$, where $n$ is the rank of $F$.

\end{lemma}
\pf 
For (1), if $x^*\in Y^\perp$ then $(x^*\otimes x) Y = 0 $ so $x^*\otimes x 
\in \calf_Y(X) $. If $x\in Y$ then $(x^*\otimes x) X \subseteq \text{span} 
\, \{x\}  \subseteq Y $ so   $x^*\otimes x \in \calf_Y(X) $.  This gives 
``$\Leftarrow$". On the other hand, if $x^* \not\in Y^\perp$ and $x\not\in 
Y$, then there is $y \in Y$ such that 
$\langle x^*, y \rangle \not= 0$, hence $(x^*\otimes x) y = \langle x^*, y 
\rangle x \not\in Y$, whence  $x^*\otimes x \not\in  \calf_Y(X) $. This 
gives ``$\Rightarrow$".

For (2), let $x_1,\dots,x_m$ be a basis for $FX\cap Y$ and extend this to 
a basis for $FX$ by adding $x_{m+1},\dots, x_n$, so that
\begin{equation*}\label{eq1}
\text{span} \ x_{m+1},\dots, x_n \cap Y = \{0\}.
\tag{*}
\end{equation*}
Write $F= \sum_{k=1}^n x_k^*\otimes x_k$.  By part (1) of this lemma, 
 $ x_k^*\otimes x_k \in 
\calf_Y(X) $ for $k\le m$. To complete the proof it is by (1) sufficient to show 
that for all $k>m$ we have  $x_k^*\in 
Y^\perp$.  If for some $k>m$ we had $x_k^*\not\in 
Y^\perp$, then, choosing $y\in Y$ with $\langle x_k^*, y\rangle \not= 0$, 
we would have by  (\ref{eq1}) that
$\sum_{j=m+1}^n \langle x_j^*, y\rangle x_j \not\in Y$.  But $\sum_{j=1}^m 
\langle x_j^*, y\rangle x_j \in Y$, so we would have $Fy \not\in Y$, a 
contradiction.
\qed

\bigskip
Theorem \ref{JointAPthm} is the main result of this note.  $N(X,Z)$ 
denotes the nuclear operators from $X$ to $Z$ and is abbreviated as $N(X)$ 
when $X=Z$. In the hypothesis we assume that the space $X$ has the AP in 
order to formulate the theorem with nuclear operators $N(X)$ rather than 
with the projective tensor product of $X^*$  with $ X$.
Given $T\in N(X)$, $\tr(T)$ is the trace of $T$, which is well-defined 
when $X$ has the AP by Grothendieck's fundamental result \cite{Gr}, 
\cite[Theorem 1.e.15]{LT}.     

\begin{theorem}\label{JointAPthm}
Suppose that $Y\subseteq X$ and $X$ has the AP. The following are 
equivalent.
\item{1.} The pair $(X,Y)$ has the AP.
\item{2.} For all $T\in N(X)$ for which $TX\subseteq Y$ and $TY=0$ we have 
$\tr(T)=0$.

\end{theorem}

\pf

(1) $\implies$ (2). Assume (2) is false and get $T\in N(X)$ so that $TY=0$ 
and $TX\subseteq Y$ but $\tr(T) = 1$.  So $T\in L(X, \tau)^*$ and $\langle 
I, T \rangle =1$. Let $F \in \calf_Y(X)$.  We want to show that $\langle 
F, T\rangle =0$, which would contradict (1).  By Lemma \ref{lemmaRank1}, 
it is enough to check that $\langle x^*\otimes x, T \rangle =0$ if either 
$x^*\in Y^\perp$ or $x \in Y$.  But $ \langle x^*\otimes x, T \rangle = 
\langle x^*, Tx\rangle$, so this is clear from the facts that  $Tx \in Y$ 
and $TY =0$.

(2) $\implies$ (1). If (1) is false and $X$ has the AP, we can separate 
$I$ from $\calf_Y(X) $ with a $\tau$ continuous linear functional on $L(X)$, 
which, since $X$ has the AP, is represented by a
nuclear operator $T$ on $X$ (see \cite{Gr}, \cite[Theorems 1.e.3, 1.e.4]{LT}). 
Then $\tr(T)= \langle I, T \rangle \not= 0$ but $\langle F, T \rangle = 0$ 
for all $F\in \calf_Y(X) $.  In particular, $\langle x^*, Tx \rangle = 
\langle x^*\otimes x, T \rangle =0$ if either $x^* \in Y^\perp$ or $x\in 
Y$.  So if $x\in X$, then for all $x^*\in Y^\perp$ we have $\langle x^*, 
Tx \rangle =0$, which is to say that $Tx \in (Y^\perp)_\perp = Y$.  So $TX 
\subseteq Y$.   If $y\in Y$, then for all $x^* \in X^*$ we have $\langle 
x^*, T y \rangle = \langle x^* \otimes y, T \rangle = 0$, which says that 
$TY=0$.
\qed

\bigskip
A sequence $Y\to X \to Z$ of Banach spaces is a short exact sequence (ses) 
when the operator $Y\to X$ is an isomorphic embedding and the operator 
$X\to Z$ is surjective and has $Y$ as its kernel.  Up to passing to 
equivalent norms, this is just saying that $Y$ is a subspace of $X$ and 
$Z$ is the quotient space $X/Y$.  A ses $Y\to X \to Z$  locally splits if 
the dual ses $Z^*\to X^* \to Y^*$ splits, which just means that $Z^*$ is a 
complemented subspace of $X^*$. This is equivalent to saying that finite 
dimensional subspaces of $Z$ uniformly lift to $X$.  The theory of ses of 
Banach spaces is presented in \cite{CG}, but much more than we use is 
contained in  \cite[Corollary 1.4]{J} and the discussion preceding that 
Corollary.  If $X\to X^{**}$ is the natural embedding, then the ses $X\to 
X^{**}\to X^{**}/X$ locally splits, but $X\to X^{**}\to X^{**}/X $  need 
not split (e.g., $X=c_0$).

\begin{prop} \label{locallySplit}
Suppose that $Y \to X \to X/Y$ is a short exact sequence that locally 
splits,  and $Y^{** *}$ and $X$ both have the AP.  Then the pair $(X,Y) $ 
has the AP and hence $X/Y$ has the AP.
\end{prop}

\pf
By Theorem \ref{JointAPthm} it is enough to show that if $T$ is a nuclear 
operator on $X$ such that $TX \subseteq Y$ and $TY = 0$, then the trace of 
$T $ is zero.   Consider $T$ as an operator into $X^{**} = Y^{\perp\perp} 
\oplus Z$ (where $Z$ is isomorphic to $(X/Y)^{**}$). This is also nuclear, 
and composing with the projection of $X^{**}$ onto $Y^{\perp\perp}(\equiv Y^{**})$ we see 
that $T$ is also nuclear when considered as an operator into $Y^{**}$. 
Since $Y^{** *}$  has the AP,  by the corrected theorem of Grothendieck 
\cite{Gr} proved by Oja and Reinov  \cite{OR}, $T$ is nuclear when 
considered as an operator from $X$ into $Y$.  Since $T$ is zero on $Y$, the 
trace of $T$ is zero.
\qed

\bigskip
Remark. The assumption on $Y$ in Proposition \ref{locallySplit} cannot be 
weakened to ``$Y^{**}$ has the AP".  (Consider a James-Lindenstrauss $Y$ 
such that  $Y^{**}$ has a basis and $Y^{**}/Y$ is a reflexive space that 
fails the AP, and let $X=Y^{**}$. It was this kind of example that led Oja 
and Reinov to the correct statement of Grothendieck's ``theorem".)

\begin{corollary} \label{Linfty}
Suppose that $Y \to X \to X/Y$ is a short exact sequence,     $X$ has the 
AP, and $Y$ is $\call_\infty$.  Then the pair $(X,Y) $ has the AP and 
hence $X/Y$ has the AP.
\end{corollary}

\pf 
The short exact sequence locally splits because $Y$ is $\call_\infty$. 
The space  $Y^{***}$ is  $\call_1$ and thus has the AP, so  the conclusion 
follows from Proposition \ref{locallySplit}.
\qed

\bigskip
In Corollary \ref{Linfty}, the roles of $X$ and $Y$ can be interchanged.

\begin{corollary} \label{Linfty2}
Suppose that $Y \to X \to X/Y$ is a short exact sequence,     $X$  is 
$\call_\infty$, and $Y$ has the AP.  Then the pair $(X,Y) $ has the AP and 
hence $X/Y$ has the AP.
\end{corollary}

\pf
By Theorem \ref{JointAPthm} it is enough to show that if $T$ is a 
nuclear operator on $X$  such that  $TX \subseteq Y$ and $TY = 0$, then the trace 
of $T $ is zero.   Just as in Proposition \ref{locallySplit}, for that it 
is enough to check that $T$ is nuclear when considered as an operator into 
$Y$.  It is, by an observation of Stegall and Retherford \cite[Theorem 
III.3]{SR}, because $X$ is $\call_\infty$.
\qed

\bigskip
Remark.  The BAP version of Corollary \ref{Linfty2} was proved in 
\cite{FJP}.

\bigskip
Next we prove  a
BAP version of Corollary \ref{Linfty}. It gives a slight improvement of 
the Castillo-Moreno \cite[Lemma 3.1]{CM}  result that $X/Y$ has the BAP if 
$X$ has the BAP and
$Y$ is $\call_\infty$. Our direct geometrical argument gives an alternate 
proof of the Castillo-Moreno result but is clumsier than the algebraic 
argument in \cite{CM} or the proof above of Corollary \ref{Linfty}.

\begin{prop} \label{Linfty3}
Suppose that $Y \to X \stackrel{Q}\to X/Y$ is a short exact sequence, 
$X$ has the BAP, and $Y$ is $\call_\infty$.  Then the pair $(X,Y) $ has 
the BAP.
\end{prop}

\pf 
Let $G$ be a finite dimensional subspace of $X$. We want to find $T\in 
\calf_Y(X)$ that is the identity on $G$ and has ``good" norm (here and in 
the following ``good" means that the norm is independent of $G$). Since 
$Y$ is $\call_\infty$, the short exact sequence locally  splits, so there 
is an operator $U: QG \to X$ with $QU$ the identity on $QG$ and $\|U\|$ 
depends only on how well the short exact sequence locally splits. So $Y + 
UQG$ is a ``good" direct sum decomposition of the space $Y + UQG$ because $U$ 
is a ``good" isomorphism on $QG$ and $QU$ is the identity on the range of 
$QG$.  By basic linear algebra, there is a finite dimensional subspace $E$ 
of $Y$ such that $G\subseteq E+UQG$.
Now $X$ has the BAP, so there is $S\in \calf(X)$ with $S$ the identity on 
$UQG$ and the norm of $S$ controlled by the BAP constant of $X$.

We next replace $S$ with an operator $S_1\in \calf(X)$ that is still the 
identity on $UQG$, has controlled norm, and is zero on $Y$ (so that  $S_1\in \calf_Y(X)$). 
 Since $G$ is arbitrary and 
$UQG=QG$, this will give the Castillo-Moreno lemma mentioned above. To get 
$S_1$, we define an operator $V\in \calf(X)$ that agrees with $S$ on $Y$ 
with $SV$ vanishing on $UQG$ so that $\|V\|$ is controlled and set $S_1:= 
S - V$.  The $\call_\infty$ structure of $Y$ is used to define $V$. 
Write $Y$ as a directed union of a net $Y_\alpha$ of subspaces of $Y$ so 
that  the   $Y_\a$ are uniformly isomorphic to  $\ell_\infty^{n_\alpha}$; 
$n_\alpha<\infty$.  Since $Y + UQG$ is a ``good" direct sum decomposition, 
the projections $P_\alpha$ from $Y_\alpha + UQG$ onto $Y_\alpha$ that are 
zero on $UQG$ have uniformly bounded norm, and, by the injective property 
of $\ell_\infty$ spaces, these projections extend to uniformly bounded 
projections (still denoted by $P_\alpha$) from $X$ onto $Y_\alpha$.  Of 
course, the net $(P_\alpha)$ converges pointwise on $Y$ to the identity on 
$Y$, so the net $SP_\alpha$ has (since $S$ has finite rank) a subnet that 
converges pointwise on $X$, necessarily to a finite rank operator $V$ that 
agrees with $S$ on $Y$ and is zero on $UQG$.
This completes the construction of $S_1:= S - V$ (and, incidentally, our 
alternate proof
for [CM, Lemma 3.1]). The remainder of the proof is very easy.  Just 
take $\alpha $ so that $E$ is a subspace of $Y_\alpha$  and define $T:= 
P_\alpha + S_1$.
\qed


\end{document}